\documentclass[12pt,noinfoline]{imsart}

\RequirePackage[OT1]{fontenc}
\RequirePackage{amsthm,amsmath,amsfonts,amssymb}
\RequirePackage[colorlinks,citecolor=blue,urlcolor=blue]{hyperref}
\usepackage{latexsym,epsf,graphicx}
\usepackage{dsfont,color,url,eurosym}
\pubyear{Date: \today}

\startlocaldefs

\newcommand{\BlackBox}{\rule{1.5ex}{1.5ex}}  

\newtheorem{theorem}{Theorem}

\newtheorem{remark}[theorem]{Remark}
\newtheorem{corollary}[theorem]{Corollary}
\newtheorem{definition}[theorem]{Definition}

\endlocaldefs

\usepackage{amsmath,amsmath,amsfonts,amssymb}
\newtheorem{assumption}[theorem]{Assumption}
\newenvironment{proofof}[1]{\begin{list}{{\bf \em Proof of #1. }}%
{\setlength{\labelsep}{0pt}\setlength{\leftmargin}{0pt}\setlength{\labelwidth}{0pt}}\item}%
{\hfill\BlackBox\end{list}}
\usepackage{dsfont} 
\usepackage{mathrsfs}
\usepackage{amsmath,amsfonts,latexsym,epsf,graphicx}

\newcommand{\bc}{\begin{center}}
\newcommand{\ec}{\end{center}}
\newcommand{\bi}{\begin{itemize} \parsep0.2em \itemsep0.2em}
\newcommand{\ei}{\end{itemize}}
\newcommand{\bnum}{\begin{enumerate} \parsep0.2em \itemsep0.2em}
\newcommand{\enum}{\end{enumerate}}
\newcommand{\be}{\begin{equation}}
\newcommand{\ee}{\end{equation}}
\newcommand{\beq}{\begin{eqnarray}}
\newcommand{\eeq}{\end{eqnarray}}
\newcommand{\beqna}{\begin{eqnarray*}}
\newcommand{\eeqna}{\end{eqnarray*}}
\newcommand{\bd}{\begin{displaymath}}
\newcommand{\ed}{\end{displaymath}}
\newcommand{\bt}{\begin{tabular}}
\newcommand{\et}{\end{tabular}}

\newcommand{\Ex}{{\mathbb{E}}}

\newcommand{\Q}{{\rm Q}}

\newcommand{\Law}[2]{\mathfrak{L}_{#1}(#2)}
\newcommand{\fL}{\mathfrak{L}}

\newcommand{\PM}  {\mathcal{M}_1}   
\newcommand{\PMZ} {\PM(\cZ,\B(\cZ))} 
\newcommand{\PMW} {\PM(\cW,\B(\cW))}

\newlength{\fixboxwidth}
\newcommand{\N}{\mathds{N}}    
\newcommand{\R}{\mathds{R}}    

\def \P {{\textrm{P}}}    
\def \Q {\textrm{Q}}    

\def \B         {{\cal B}}                

\def \lb        { \lambda }

\def \ve        { \varepsilon }

\def \s         { \sigma }

\def \et        { \tilde{\eta}}

\newcommand{\snorm}[1] {\Vert #1 \Vert}
\newcommand{\BLnorm}[1]{\snorm{#1}_{\mathrm{BL}}}

\newcommand{\inorm}[1]{\snorm{#1}_{\infty}}

\newcommand{\sLp}[1]{\mbox{${\cal L}_p(\mu)$}}

\newcommand{\cE}    {\mathcal{E}} 
 
\newcommand{\cP}    {\mathcal{P}} 
\newcommand{\cS}    {\mathcal{S}} 
\newcommand{\cX}    {\mathcal{X}}  
\newcommand{\cY}    {\mathcal{Y}}  
\newcommand{\cW}    {\mathcal{W}}  
\newcommand{\cZ}    {\mathcal{Z}} 

\newcommand{\hnorm}[1]{\left\Vert #1 \right\Vert_{H}}
\newcommand{\hhnorm}[1]{\left\Vert #1 \right\Vert_{H}^2}

\def \P           { \mathrm{P} }   
\def \Q           { \mathrm{Q} } 

\newcommand{\D} {\mathrm{D}}

\newcommand{\Ls}{L^{\star}}

\newcommand{\cXY}{\cX\times\cY}

\newcommand{\cXYR}{\cX\times\cY\times\R}

\newcommand{\dBL}{{{d}_{\mathrm{BL}}}}

\newcommand{\cA}{\mathcal{A}}

\begin{document}

\begin{frontmatter}
\title{On the Stability of Bootstrap Estimators}
\runtitle{On the Stability of Bootstrap Estimators}

\begin{aug}
\author{\fnms{A.} \snm{Christmann}
\ead[label=e1]{andreas.christmann@uni-bayreuth.de}
}
\and
\author{\fnms{M.} \snm{ Salib\'{i}an-Barrera}
\ead[label=e2]{matias@stat.ubc.ca}}
\and
\author{\fnms{S.} \snm{Van Aelst}
\ead[label=e3]{Stefan.VanAelst@UGent.be}}

\address[a]{University of Bayreuth, Department of Mathematics, Bayreuth, GERMANY.\\
\printead{e1}}

\address[b]{University of British Columbia, Department of Statistics, Vancouver, CANADA.\\ \printead{e2}}

\address[c]{University of Ghent, Department of Applied Mathematics and Computer Science, Ghent, BELGIUM.\\ \printead{e3}}

\runauthor{Christmann, Salib\'{i}an-Barrera, Van Aelst}

\affiliation{University of Bayreuth, University of British Columbia, University of Ghent}
\end{aug}

\begin{abstract}
It is shown that bootstrap approximations of an estimator which is based on a continuous operator from the set of Borel probability measures defined on a compact metric space into a complete separable metric space is stable in the sense of qualitative robustness.
Support vector machines based on shifted loss functions are treated as special cases.
\end{abstract}


\end{frontmatter}


\section{Introduction}
  The finite sample distribution of many 
nonparametric methods from statistical learning theory
is unknown because the distribution $\P$ from which the 
data were generated is unknown and because there are often
only asymptotical results on the behaviour of such methods known.

The goal of this paper is to show that bootstrap approximations of an estimator which is based on a continuous operator from the set of Borel probability distributions defined on a compact metric space into a complete separable metric space is stable in the sense of qualitative robustness.
As a special case it is shown that bootstrap approximations for the
support vector machine (SVM) are stable, both for the risk functional and for the SVM operator itself.
The results can be interpreted as generalizations of theorems derived by 
\cite{CuevasRomo1993}.

The rest of the paper has the following structure.
Section 2 gives the general result and Section 3 contains the results for SVMs.
All proofs are given in the appendix.

\section{On Qualitative Robustness of Bootstrap Estimators}\label{sec:qualirob}

If not otherwise mentioned, we will use the Borel $\s$-algebra $\B(A)$ on 
a set $A$ and denote the Borel $\s$-algebra on $\R$ by $\B$.

\begin{assumption} \label{AC.genass}
Let $(\Omega,\cA,\mu)$ be a probability space, where $\mu$ is unknown,
$(\cZ,d_\cZ)$ be a compact metric space, and 
$\B(\cZ)$ be the Borel $\s$-algebra on $\cZ$.
Denote the set of all Borel probability measures on
$(\cZ,\B(\cZ))$ by $\PM(\cZ,\B(\cZ))$. On $\PM(\cZ,\B(\cZ))$ we use the Borel $\s$-algebra $\B(\PM(\cZ,\B(\cZ)))$ and the bounded Lipschitz metric $\dBL$, see {(\ref{AC.Lemma1starC})}. 
Let $S$ be a statistical operator defined on 
$\PM(\cZ,\B(\cZ))$ 
with values in a complete, separable metric space $(\cW,d_\cW)$ enclipped with its Borel $\s$-algebra $\B(\cW)$.
Let $Z,Z_n: (\Omega,\cA,\mu) \to (\cZ,\B(\cZ))$, $n\in\N$, 
be independent and identically distributed random variables and denote the
image measure by $\P:=Z\circ\mu$.
Let $S_n(Z_1,\ldots,Z_n)$ be a statistic with values in $(\cW, \B(\cW))$.
Denote the empirical measure of $(Z_1,\ldots,Z_n)$ by
$\P_n:=\frac{1}{n}\sum_{i=1}^n \delta_{Z_i}$. 
The statistic $S_n$ is defined via the operator
$$S: (\PM(\cZ,\B(\cZ)), \B(\PM(\cZ,\B(\cZ))) \to (\cW,\B(\cW))$$ 
where $S(\P_n)=S_n(Z_1,\ldots,Z_n)$.
Denote the distribution of $S_n(Z_1,\ldots,Z_n)$ when 
$Z_i \stackrel{i.i.d.}{\sim} \P$ by
$\Law{n}{S;\P}:=\fL(S_n(Z_1,\ldots,Z_n))$.
Accordingly, we denote the distribution of $S_n(Z_1,\ldots,Z_n)$ when 
$Z_i \stackrel{i.i.d.}{\sim} \P_n$ by
$\Law{n}{S; \P_n}$. 
\end{assumption}

Efron \cite{Efron1979,Efron1982} proposed the bootstrap, whose main idea is 
to approximate the unknown distribution $\Law{n}{S;\P}$ 
by $\Law{n}{S;\P_n}$. 
Note that these bootstrap approximations $\Law{n}{S; \P_n}$ are
(probability measure-valued) random variables with values in 
$\PMW$.

Following \cite{CuevasRomo1993} we call a sequence of bootstrap approximations 
$\Law{n}{S;\P_n}$ \emph{qualitatively robust} at $\P\in\PM(\cZ,\B(\cZ))$ 
if the sequence of transformations 
\be
 g_n: \PMZ  \to \PMW,
 \quad
 g_n(\Q) = \fL(\Law{n}{S;\Q_n}),
 \qquad n\in\N,
\ee
is asymptotically equicontinuous at $\P\in\PMZ$, i.e. if
\be
\begin{split}
\forall \, \ve>0 ~\exists\, \delta>0 ~\exists \, n_0\in\N :
\hfill \qquad \qquad \qquad \qquad \qquad \qquad \qquad \qquad ~\\
\dBL(\Q,\P) < \delta ~~ \Rightarrow ~~
  \sup_{n\ge n_0} \dBL\bigl(\fL(\Law{n}{S;\Q_n}), \fL(\Law{n}{S;\P_n})\bigr) < \ve.
\end{split}
\ee
Following \cite{CuevasRomo1993} again, we call a sequence of statistics $(S_n)_{n\in\N}$ \emph{uniformly qualitatively robust in a neighborhood $\mathcal{U}(\P_0)$ of $\P_0\in\PMZ$} if 
\be 
\begin{split} 
\exists \, n_0\in\N ~\forall \, \ve>0 ~\forall \, n\ge n_0
~\exists\, \delta>0 
~\forall\, \P \in \mathcal{U}(\P_0): \quad~ \\
  \dBL(\Q,\P) < \delta \quad \Rightarrow \quad
  \dBL(\Law{n}{S;\Q}, \Law{n}{S;\P}) < \ve.
\end{split}
\ee

The following two results and Theorem
\ref{AC.TheoremQualiBootSVM} in the next section 
are the main results of this paper.

\begin{theorem}\label{AC.Theorem2star}
If Assumption \ref{AC.genass} is valid and  
if $S$ is uniformly continuous in a neighborhood 
$\mathcal{U}(\P_0)$ of $\P_0\in\PM(\cZ,\B(\cZ))$, 
then $(S_n(Z_1,\ldots,Z_n))_{n\in\N}$ is uniformly qualitatively robust in $\mathcal{U}(\P_0)$.
\end{theorem}

\begin{theorem}\label{AC.Theorem3star}
If Assumption \ref{AC.genass} is valid and  
if $(S_n(Z_1,\ldots,Z_n))_{n\in\N}$ is uniformly qualitatively robust in a neighborhood 
$\mathcal{U}(\P_0)$ of 
$\P_0\in\PM(\cZ,\B(\cZ))$, then the sequence $\Law{n}{S;\P_n}$
of bootstrap approximations of  $\Law{n}{S;\P}$ is qualitatively robust for $\P_0$. 
\end{theorem}

As an immediate consequence from both theorems given above we obtain

\begin{corollary}\label{AC.Corollary3star}
If Assumption \ref{AC.genass} is valid and if $S$ is a continuous operator, then
the sequence $\Law{n}{S;\P_n}$ of bootstrap approximations of $\Law{n}{S;\P}$ is qualitatively robust for all $\P \in \PM(\cZ,\B(\cZ))$. 
\end{corollary}

\begin{remark}
The Theorems \ref{AC.Theorem2star} and \ref{AC.Theorem3star} can be considered as a generalization of \cite[Thm.\,2, Thm. 3]{CuevasRomo1993}, who considered the case  $\cW:=A\subset\R$ being a finite interval and $\cZ:=\R$-valued random variables $Z_1,\ldots,Z_n$. 
In our case, the statistics $S_n(Z_1,\ldots,Z_n)$ are $\cW$-valued statistics, where $\cW$ is a complete separable metric space and its dimension can be infinite.
\end{remark}


\section{On Qualitative Robustness of Bootstrap SVMs}\label{sec:qualirobsvm}

In this section we will apply the previous results to support vector machines which belong to the modern class of statistical machine learning methods. I.e., we will consider the special case that $\cW$ is a reproducing kernel Hilbert space $H$ used by a support vector machine (SVM). Note that $H$ typically has an infinite dimension, which is true, e.g., if the popular Gaussian RBF kernel 
$k:\cX\times\cX\to\R$, $k(x,x'):=\exp(-\gamma \snorm{x-x'}_2^2)$ for $\gamma>0$) is used.

To state our result on the stability of bootstrap SVMs in Theorem \ref{AC.TheoremQualiBootSVM} below, we need the following assumptions on the loss function and the kernel. 

\begin{assumption}\label{AC.SVMass}
Let $\cZ=\cXY$ be a compact metric space with metric
$d_\cZ$, where $\cY\subset \R$ is closed.
Let $L:\cX\times\cY\times\R\to[0,\infty)$ be a loss function such that $L$ is continuous and convex with respect to its third argument and that $L$ is 
uniformly Lipschitz continuous with respect to its third argument with uniform Lipschitz constant $|L|_1>0$, i.e.
$|L|_1$ is the smallest constant $c$ such that
$\sup_{(x,y)\in \cXY}|L(x,y,t)-L(x,y,t')| \le c |t-t'|$ for all $t,t'\in\R$. Denote the shifted loss function by
$\Ls(x,y,t):=L(x,y,t)-L(x,y,0)$, $(x,y,t)\in \cXYR$.
Let $k:\cX\times\cX\to\R$ be a continuous kernel with reproducing kernel Hilbert space $H$ and assume that $k$ is bounded by $\inorm{k}:=(\sup_{x\in\cX} k(x,x))^{1/2} \in (0,\infty)$.
Let $\lb\in(0,\infty)$. 
\end{assumption}

These assumptions can be considered as standard assumptions for stable SVMs, see, e.g., \cite{ChristmannSteinwart2007a} and \cite[Chap.\,10]{SteinwartChristmann2008a}, .

In this paper the RKHS $H$, the penalyzing constant $\lb$, and the loss function $L$ and thus the shifted loss function $\Ls$ are fixed. Therefore, we write in the next definition just $S$ and $R$ instead of 
$S_{\Ls, H,\lb}$ and $R_{\Ls, H, \lb}$ to shorten the notation.

\begin{definition}\label{AC.SVMdef}
The \textbf{SVM operator} 
$S: \PM(\cZ,\B(\cZ)) \to H$ is defined by
\be \label{SVMoperatoOn r}
S(\P)
:= 
f_{\Ls,\P,\lb}
:=\arg \min_{f\in H} \Ex_\P \Ls(X,Y,f(X)) + \lb \hhnorm{f}.
\ee 
The \textbf{SVM risk functional} 
$R: \PM(\cZ,\B(\cZ)) \to \R$ is defined by
\be \label{SVMfunctional}
R(\P)
:= 
\Ex_\P \Ls(X,Y,S(\P)(X))
=
\Ex_\P \Ls(X,Y,f_{\Ls,\P,\lb}(X)).
\ee 
\end{definition}

If Assumption \ref{AC.SVMass} is valid, then $S$ is 
well-defined because $S(\P) \in H$ exists and is unique, 
$R$ is well-defined because $R(\P)\in \R$ exists and is unique, and it holds, for all $\P\in\PM(\cXY)$,
\be \label{AC.SVMsbounds}
 \inorm{S(\P)} \le \frac{1}{\lb}  |L|_1 \, \inorm{k}^2 < \infty
 \quad \mathrm{and} \quad
 |R(\P)| \le \frac{1}{\lb}  |L|_1^2 \, \inorm{k}^2 < \infty \, ,
\ee
see 
\cite[Thm\,5, Thm.\,6, (17),(18)]{ChristmannVanMessemSteinwart2009}.

\begin{theorem}\label{AC.TheoremQualiBootSVM}
If the general Assumption \ref{AC.genass} and 
the Assumption \ref{AC.SVMass} are valid, then the SVM operator $S$ and the SVM risk functional $R$
fulfill:
\bnum
\item[(i)] The sequence $\Law{n}{S;\P_n}$ of bootstrap SVM estimators of $\Law{n}{S;\P}$ is qualitatively robust for all $\P\in\PM(\cZ,\B(\cZ))$.
\item[(ii)] The sequence $\Law{n}{R;\P_n}$ of bootstrap SVM risk estimators of $\Law{n}{R;\P}$ is qualitatively robust for all $\P\in\PM(\cZ,\B(\cZ))$.
\enum
\end{theorem}

\section{Proofs}
  \subsection{Proofs of the results in Section \ref{sec:qualirob}}

For the proofs we need Theorem \ref{AC.Lemma1star} and Theorem \ref{ThmStrassendBL}, see below. 
To state Theorem \ref{AC.Lemma1star} on uniform Glivenko-Cantelli classes, we need the following notation.
For any metric space $(\cS,d)$ and real-valued function 
$f:\cS\to\R$, we denote the bounded Lipschitz norm of $f$ by
\be \label{fBL}
\BLnorm{f} 
:=
\sup_{x\in\mathcal{S}} |f(x)| + 
\sup_{x,y\in\mathcal{S}, x\ne y} \frac{|f(x)-f(y)|}{d(x,y)} \,.
\ee
Let $\tilde{F}$ be a set of measurable functions from $(\cS,\B(\cS)) \to (\R,\B)$.
For any function $G:\tilde{F}\to\R$ (such as a signed measure) define
\be \label{tildeF}
\snorm{G}_{\tilde{F}} 
:=
\sup \{|G(f)|: f\in\tilde{F}\}.
\ee

\begin{theorem}\label{AC.Lemma1star}
\cite[Prop. 12]{DudleyGineZinn1991}
For any separable metric space $(\mathcal{S},d)$ 
and $M \in (0, \infty)$, 
\be \label{AC.Lemma1starA}
\tilde{\mathcal{F}}_M
:=
\{f:(\cS,\B(\cS)) \to (\R,\B); \snorm{f}_{BL} \le M\}
\ee
is a universal Glivenko-Cantelli class. It is a uniform  Glivenko-Cantelli class, 
i.e., for all $\ve > 0$, 
\be \label{AC.Lemma1starB}
\lim_{n\to\infty} \sup_{\nu\in\PM(\cS,\B(\cS))} 
{\rm Pr}^*\Bigl(\sup_{m\ge n} \snorm{\nu_m-\nu}_{\tilde{\mathcal{F}}_M} > \ve\Bigr) =0,  
\ee
if and only if $(\mathcal{S},d)$ is totally bounded.
Here, ${\rm Pr}^*$ denotes the outer probability.
\end{theorem}

Note that the term $\snorm{\nu_m-\nu}_{\tilde{\mathcal{F}}_M}$ in
{(\ref{AC.Lemma1starB})} equals the \emph{bounded Lipschitz metric} $\dBL$ of the probability measures $\nu_m$ and 
$\nu$ if $M=1$, i.e.
\be \label{AC.Lemma1starC}
\snorm{\nu_m-\nu}_{\tilde{\mathcal{F}}_1}
=
\sup_{f\in\tilde{F}_1} |(\nu_m-\nu)(f)|
=
\sup_{f; \BLnorm{f}\le 1} \Bigl| \int \!f\,d\nu_m - \int \!f \,d\nu \Bigr|
=:
\dBL(\nu_m,\nu),
\ee
see \cite[p.\,394]{Dudley2002}.
Hence, Theorem \ref{AC.Lemma1star} can be interpreted as a generalization of 
\cite[Lemma 1, p. 186]{CuevasRomo1993}, which says that
if $A\subset\R$ is a finite interval, then $\dBL(\P_m,\P)$ converges almost surely to $0$ uniformly in $\P\in\PM(A,\B(A))$.
For various characterizations of Glivenko-Cantelli classes, we refer to
\cite[Thm. 22]{Talagrand1987} and \cite{Dudley1999}.

We next list the other main result we need for the proof of 
Theorem \ref{AC.TheoremQualiBootSVM}. This result is an analogon of the famous Strassen theorem for the bounded Lip\-schitz metric $\dBL$ instead of the Prohorov metric.

\begin{theorem}\label{ThmStrassendBL} \cite[Thm. 4.2, p. 30]{Huber1981}
Let $\cZ$ be a Polish space with topology $\tau_\cZ$.
Let $\dBL$ be the bounded Lipschitz metric defined on the 
set $\PM(\cZ,\B(\cZ))$ of all Borel probability measures on $\cZ$.
Then the following two statements are equivalent:
\bnum
\item[(i)] There are random variables $\xi_1$ with distribution $\nu_1$
and $\xi_2$ with distribution $\nu_2$ such that
$\Ex[\dBL(\xi_1, \xi_2)] \le \ve$.
\item[(ii)] $\dBL(\nu_1, \nu_2) \le \ve$.
\enum
\end{theorem}

\begin{proofof}{Theorem \ref{AC.Theorem2star}}
We closely follow the proof by 
\cite[Thm. 2]{CuevasRomo1993}. However, we use  
Theorem \ref{AC.Lemma1star} instead of their Lemma 1
and we use \cite[Lem. 1]{Cuevas1988} instead of 
\cite[Lem. 1]{Hampel1971}.

Let $\cP_n\subset \PMZ$ be the set of empirical distributions of order $n\in\N$,
i.e. 
\be
\cP_n := 
\Bigl\{
\P_n \in \PMZ; \exists\, (z_1,\ldots,z_n)\in\cZ^n \mathrm{~such~that~}
\P_n=\frac{1}{n}\sum_{i=1}^n \delta_{z_i}
\Bigr\}\,,
\ee
and let $\cE_n \subset \cP_n$. If misunderstandings are unlikely, we identify $\cE_n$ with the set $\{z_1,\ldots,z_n\}$ of atoms.

It is enough to show that
\be
\forall\,\ve>0 ~\exists\,\delta>0 
~\forall\,\P\in\mathcal{U}(\P_0)
~\exists\,\mathrm{sequence~} (\mathcal{E}_n)_{n\in\N}\subset\mathcal{P}_n
\ee
such that $\P^n(\mathcal{E}_n) > 1-\ve$
and
for all $\Q_n\in\mathcal{E}_n$ and for all $\tilde{\Q}_n\in\mathcal{P}_n$ we have
\be
\dBL(\Q_n,\tilde{\Q}_n) < \delta 
\quad \Rightarrow \quad
d_\cW(S(\Q_n),S(\tilde{\Q}_n)) < \ve. 
\ee
From this we obtain that $(S_n)_{n\in\N}$ is uniformly
qualitatively robust by \cite[Lem.\,1]{Cuevas1988}.

Let $\ve>0$. Since the operator $S$ is uniformly
continuous in $\mathcal{U}(\P_0)$ we obtain
\be \label{AC.dWsmall}
\exists\,\delta_0>0 
~\forall\,\P\in\mathcal{U}(\P_0):
~~
\dBL(\P,\Q)<\delta_0 
\quad  \Rightarrow \quad
d_\cW (S(\P),S(\Q))< \ve / 2 \, . 
\ee
Hence by Theorem \ref{AC.Lemma1star} for the special case $M=1$ and by {(\ref{AC.Lemma1starC})}, 
we get
\be
\exists\,n_0\in\N: 
~~
\sup_{\P\in\mathcal{U}(\P_0)} 
{\rm Pr}^*\Bigl(\sup_{n\ge n_0} \dBL(\P_n,\P) < \delta_0 \Bigr) > 1-\ve.
\ee
For $n \ge n_0$ and $\P\in\mathcal{U}(\P_0)$, define
\be
\mathcal{E}_{n,\P}
:=
\{ \Q_n\in\mathcal{P}_n: \dBL(\Q_n,\P) < \delta_0 /2 \} \,.
\ee
It follows, that
$\P^n(\mathcal{E}_{n,\P})>1-\ve$ together with
$\Q_n\in\mathcal{E}_{n,\P}$ and 
$\dBL(\Q_n,\tilde{\Q}_n) <  \delta_0 / 2$
implies that 
$$
  \dBL(\Q_n,\P) <  \delta_0 / 2
  \quad \mathrm{and} \quad
  \dBL(\tilde{\Q}_n,\P) <  \delta_0 \, .
$$
The triangle inequality thus yields due to 
{(\ref{AC.dWsmall})}
\be
d_\cW(S(\Q_n),S(\tilde{\Q}_n))
\le
d_\cW(S(\Q_n),S(\P))
+
d_\cW(S(\P),S(\tilde{\Q}_n))
< 
\ve,
\ee
from which the assertion follows.
\end{proofof}

\begin{proofof}{Theorem \ref{AC.Theorem3star}}
The proof mimics the proof of \cite[Thm. 3]{CuevasRomo1993},
but uses Theorem \ref{AC.Lemma1star} instead of 
\cite[Lem.\,1]{CuevasRomo1993}.

Fix $\P_0\in\PMZ$ and $\ve>0$.
By the uniform qualitative robustness of $(S_n)_{n\in\N}$ in 
$\mathcal{U}(\P_0)$, there exists $n\in\N$ such that for all $\ve>0$ there exists $\delta>0$ such that
\be \label{AC.Theorem3star.p1}
\dBL(\Q,\P)<\delta 
\quad  \Rightarrow \quad
\sup_{m\ge n} \sup_{\P\in\mathcal{U}(\P_0)}
   \dBL(\Law{m}{S;\Q}, \Law{m}{S;\P}) < \ve. 
\ee
Define $\delta_1:=\delta/2$.
Due to Theorem \ref{AC.Lemma1star} for the special case $M=1$ and by {(\ref{AC.Lemma1starC})}, 
we have, for all $\ve> 0$,
\be
\lim_{n\to\infty} 
\sup_{\P\in\PM(\cZ,\B(\cZ))} 
{\rm Pr}^*\Bigl(\sup_{m \ge n} \dBL(\P_m,\P) > \ve \Bigr) = 0.
\ee
Hence {(\ref{AC.Theorem3star.p1})} and Varadarajan's theorem on the almost sure convergence of empirical measures to a Borel probability measure defined on a separable metric space, see e.g. 
\cite[Thm. 11.4.1, p.\,399]{Dudley2002},  yields for the empirical distributions $\Q_n$ from $\Q$ and $\P_{0,n}$ from $\P_0$ that,  
\be \label{AC.TheGlivenkoorem3star.p2}
\exists\, n_1 > n
~\forall\, n\ge n_1:~ 
\dBL(\Q,\P_0) < \delta_1
\quad \Rightarrow \quad 
\dBL(\Q_n, \P_{0,n}) < \delta
\mathrm{~~almost~surely.}
\ee
It follows from the uniform qualitative robustness of $(S_n)_{n\in\N}$,
see {(\ref{AC.Theorem3star.p1})}, that
\be \label{AC.Theorem3star.p3}
\begin{split} 
\exists \, n_1\in\N
~\forall \, \ve>0 
~\forall \, n\ge n_1
~\exists\, \delta>0 
~\forall\, \P \in \mathcal{U}(\P_0): \qquad\qquad\qquad\qquad\quad~ \\
\dBL(\Q,\P) < \delta 
\quad \Rightarrow \quad
\dBL(\Law{n}{S;\Q_n}, \Law{n}{S;\P_{0,n}}) < \ve 
\mathrm{~~almost~surely.}
\end{split}
\ee
For notational convenience, we write for the sequences of 
bootstrap estimators 
\be 
\xi_{1,n}:=\Law{n}{S;\Q_n}, 
\qquad
\xi_{2,n}:=\Law{n}{S;\P_{0,n}}, 
\qquad n\in\N.
\ee
Note that $\xi_{1,n}$ and  $\xi_{2,n}$ are (measure-valued) random variables with values in the set $\PMW$.
We denote the distribution of $\xi_{j,n}$ by $\mu_{j,n}$ for 
$j\in\{1,2\}$ and $n\in\N$.
Hence {(\ref{AC.Theorem3star.p3})} yields
\be
\dBL(\xi_{1,n}, \xi_{2,n}) < \ve 
\mathrm{~~almost~surely~for~all~} n \ge n_1
\ee
and it follows
\be
\Ex[\dBL(\xi_{1,n}, \xi_{2,n})] \le \ve, 
\qquad \forall\,n \ge n_1.
\ee
Now an application of an analogon of Strassen's theorem,
see Theorem \ref{ThmStrassendBL}, yields
\be 
\sup_{n\ge n_1} \dBL(\mathfrak{L}(\xi_{1,n}),\mathfrak{L}(\xi_{2,n})) \le \ve
\qquad
\forall\,n\ge n_1,
\ee
which completes the proof, because 
\be
\mathfrak{L}(\xi_{1,n}) = \mathfrak{L}(\Law{n}{S;\Q_n})
\quad \mathrm{and} \quad
\mathfrak{L}(\xi_{2,n}) = \mathfrak{L}(\Law{n}{S;\P_{0,n}}).
\ee
\end{proofof}

\subsection{Proofs of the results in Section \ref{sec:qualirobsvm}}
\begin{proofof}{Theorem \ref{AC.TheoremQualiBootSVM}}
\emph{Proof of part (i).}
By assumption, $(\cZ,d_\cZ)$ is a compact metric space, where $\cZ=\cX\times\cY$. Let $\B(\cZ)$ be the
Borel $\s$-algebra on $\cZ$.
It is well-known that the bounded Lipschitz metric $\dBL$ metrizes the weak topology on the space $\PM(\cZ,\B(\cZ))$, 
see \cite[Thm. 11.3.3]{Dudley2002}, and that $(\PM(\cZ,\B(\cZ)), \dBL)$ is a compact metric space 
if and only if $(\cZ,d_\cZ)$ is a compact metric space,
see \cite[p.\,45, Thm. 6.4]{Parthasarathy1967}.  
From the compactness of $(\PM(\cZ,\B(\cZ)),\dBL)$, it
of course follows that this metric space is separable and totally bounded, see \cite[Thm. 1.4.26]{DenkowskiEtAl2003}.

Under the assumptions of the theorem we have, for all fixed $\lb\in(0,\infty)$, that the SVM operator
$S: \PM(\cZ,\B(\cZ)) \to H$, $S(\P)=f_{\Ls,\P,\lb}$, is
well-defined because it exists and is unique,
see \cite[Thm.\,5, Thm.\,6]{ChristmannVanMessemSteinwart2009} and is continuous with respect to the combination of the weak topology on $\PM(\cZ,\B(\cZ))$ and 
the norm topology on $H$, see 
\cite[Thm. 3.3, Cor. 3.4]{HableChristmann2011}. 
There it was also shown that the operator
$\tilde{S}: \PM(\cZ,\B(\cZ)) \to \mathcal{C}_b(\cZ)$, $\P \mapsto f_{\Ls,\P,\lb}$, 
is continuous with respect to the combination of weak topology on $\PM(\cZ,\B(\cZ))$ and the norm topology on $\mathcal{C}_b(\cZ)$.
Because $(\PM(\cZ,\B(\cZ)),\dBL)$ is a \emph{compact} metric space, the operators $S$ and $\tilde{S}$ are therefore even \emph{uniformly continuous} on the whole space $\PMZ$ with respect to  the mentioned topologies, see \cite[Prop. 1.5.9]{DenkowskiEtAl2003}. 

Because the reproducing kernel Hilbert space $\cW:=H$ is a Hilbert space, $H$ is complete.
Furthermore, because the input space $\cX$ is separable and the kernel $k$ is continuous, the RKHS $H$ is also separable, see 
\cite[Lem. 4.33]{SteinwartChristmann2008a}. 
Therefore, Theorem \ref{AC.Theorem2star} yields that the sequence of $H$-valued statistics 
\be 
 S_n((X_1,Y_1),\ldots,(X_n,Y_n))
 = \arg \min_{f\in H} \frac{1}{n} \sum_{i=1}^n
 \Ls(X_i,Y_i,f(X_i)) + \lb \hhnorm{f}, ~n\in\N,
\ee
is uniformly qualitatively robust in a neighborhood $\mathcal{U}(\P_0)$ for every probability measure 
$\P_0 \in \PM(\cZ)$. 
Now we apply Theorem \ref{AC.Theorem3star}, which yields that the sequence $(\Law{n}{S;\P_n})_{n\in\N}$ of bootstrap SVM estimators of $\Law{n}{S;\P}$ is qualitatively robust for \emph{all} $\P_0\in\PM(\cZ,\B(\cZ))$, which gives the 
first assertion of the theorem.

\emph{Proof of part (ii).}
The proof consists of two steps. In Step 1 the continuity of the SVM risk functional $R$ will be shown. 
In Step 2, the Theorems \ref{AC.Theorem2star} and  \ref{AC.Theorem3star} will be used to show that the sequence $(\Law{n}{R;\P_n})_{n\in\N}$, $n\in\N$, of bootstrap SVM risk estimators is qualitatively robust.

\emph{Step 1.}
We will first show that the SVM risk functional 
$R:\PMZ\to\R$ is continuous 
with respect to the combination of the weak topology on $\PM(\cZ,\B(\cZ))$ and the standard topology on $\R$.

As mentioned in part \emph{(i)}, the assumption that $(\cZ,d_\cZ)$ is a compact metric space implies that $(\PM(\cZ,\B(\cZ)), \dBL)$ is a compact metric space and hence this space is separable and totally bounded.

Under the assumptions of the theorem, the SVM operator
$S: \PM(\cZ,\B(\cZ)) \to H$, $S(\P)=f_{\Ls,\P,\lb}$, is well-defined because $S(\P)$ exists and is unique for all $\P\in\PMZ$ and for all $\lb\in(0,\infty)$, 
see \cite[Thm.\,5, Thm.\,6]{ChristmannVanMessemSteinwart2009}. Furthermore, $S$ is continuous with respect to the combination of the weak topology on $\PM(\cZ,\B(\cZ))$ and the norm topology on $H$, see \cite[Thm. 3.3]{HableChristmann2011}. 
Hence the function 
\be
g_\P: \cX\times\cY \to \R,
\quad 
g_\P(x,y)
:= \Ls\bigl(x,y,S(\P)(x)\bigr) 
= \Ls\bigl(x,y,f_{\Ls,\P,\lb}(x)\bigr)
\ee
is well-defined.
Because the kernel $k$ is bounded and continuous, all functions $f\in H$, and hence in particular 
$S(\P)=f_{\Ls,\P,\lb}\in H$, are continuous, see e.g.
\cite[Lem.\,4.28, Lem.\,4.29]{SteinwartChristmann2008a}.
Hence the function $g_\P$ is \emph{continuous} (with respect to $(x,y)$), because the loss function $L$ and hence the shifted loss function $\Ls(x,y,t)=L(x,y,t)-L(x,y,0)$, $(x,y,t)\in\cXYR$, are continuous.
Furthermore, the function $g_\P$ is \emph{bounded}, because $(\cZ,d_\cZ)$ with $\cZ:=\cX\times\cY$ is by assumption a compact metric space, the Lipschitz continuous loss function $L$ maps from $\cXYR$ to $[0,\infty)$, and
$\inorm{S(\P)} \le \frac{1}{\lb}  |L|_1 \, \inorm{k}^2 < \infty$,
see \cite[p.\,314, (17)]{ChristmannVanMessemSteinwart2009}.
Hence $g_\P\in \mathcal{C}_b(\cZ,\R)$.
Because the bounded Lipschitz metric $\dBL$ metrizes the weak topology on 
$\PMZ$, it follows that 
\be \label{AC:gPconvergence}
\forall\,\ve_1>0
~\exists\, \delta_1>0:~~~
\dBL(\Q,\P) < \delta_1
\quad \Longrightarrow \quad
\Bigl| \int g_\P \,d\Q - \int g_\P \,d\P \Bigr| 
< \ve_1 \, .
\ee

Recall that $S:\PMZ\to H$ is continuous with respect to the combination of the weak topology on $\PM(\cZ,\B(\cZ))$ and 
the norm topology on $H$, see \cite[Thm. 3.3]{HableChristmann2011}. 
Hence
\be \label{AC:Sconvergence}
\forall\,\ve_2>0
~\exists\, \delta_2>0:~~~
\dBL(\Q,\P) < \delta_2
\quad \Longrightarrow \quad
\hnorm{S(\Q)-S(\P)} < \ve_2 \, .
\ee

Fix $\ve > 0$.
Define 
$$ 
\ve_1:=\frac{\ve}{3}
\quad \mathrm{and} \quad 
\ve_2:=\frac{\ve}{3 |L|_1 \inorm{k}} ~~.
$$
Using the triangle inequality in {(\ref{AC.Tcontinuous1})},
the definition of the shifted loss function $\Ls$ in {(\ref{AC.Tcontinuous2a})},
the definition of the function $g_\P$ in 
{(\ref{AC.Tcontinuous2b})},
the Lipschitz continuity of $L$ in {(\ref{AC.Tcontinuous3})},
and the well-known formula
\be \label{AC:inequality1}
\inorm{f} \le \inorm{k} \hnorm{f}, \qquad f\in H,
\ee 
see e.g. \cite[p.\,124]{SteinwartChristmann2008a}  we obtain that
$\dBL(\Q,\P) < \delta_2$ implies
\begin{eqnarray}
 & & |R(\Q)-R(\P)|  \nonumber\\
 & \!\!=\!\! &
 \Bigl| \int \Ls(x,y,S(\Q)(x))\,d\Q(x,y) - \int \Ls(x,y,S(\P)(x))\,d\P(x,y) \Bigr| 
 \nonumber \\
 & \!\!\le\!\! & 
  \Bigl| \int \!\Ls(x,y,S(\Q)(x))\,d\Q(x,y) 
  \!-\! \int \!\Ls(x,y,S(\P)(x))\,d\Q(x,y) \Bigr| \label{AC.Tcontinuous1} \\
 & & + \Bigl| \int \!\Ls(x,y,S(\P)(x))\,d\Q(x,y) \!-\! 
 \int \!\Ls(x,y,S(\P)(x))\,d\P(x,y) \Bigr|  \nonumber \\
 & \!\!\le\!\!  & 
  \int | L(x,y,S(\Q)(x))\!-\!L(x,y,S(\P)(x))| \,d\Q(x,y) \label{AC.Tcontinuous2a} \\
 & &  + ~ \Bigl| \int g_\P \,d\Q - \int g_\P \,d\P \Bigr|   \label{AC.Tcontinuous2b} \\ 
 & \!\!\stackrel{{\scriptsize{(\ref{AC:gPconvergence})}}}{\le}\!\! & 
  |L|_1 \, \inorm{S(\Q)-S(\P)} + \ve_1 \label{AC.Tcontinuous3} \\
 & \!\!\stackrel{{\scriptsize{(\ref{AC:inequality1})}}}{\le}\!\!  & 
  |L|_1 \, \inorm{k}\, \hnorm{S(\Q)-S(\P)} + \ve_1 \\
 & \!\!\stackrel{{\scriptsize{(\ref{AC:Sconvergence})}}}{\le}\!\! & 
 |L|_1 \, \inorm{k}\, \ve_2 + \ve_1 = \frac{2}{3} \, \ve. 
\end{eqnarray}

Hence, $R$ is continuous with respect to the combination of the weak topology on $\PM(\cZ,\B(\cZ))$ and the standard topology on $(\R,\B)$.

\emph{Step 2.}
Because $(\PM(\cZ,\B(\cZ)),\dBL)$ is a \emph{compact} metric space and the risk functional
$R:\PMZ\to\R$ is continuous, $R$ is even \emph{uniformly continuous} with respect to the mentioned topologies, see \cite[Prop. 1.5.9]{DenkowskiEtAl2003}. 
Obviously $(\cW,d_\cW):=(\R,|\cdot|)$ is a complete separable metric space.
Therefore, Theorem \ref{AC.Theorem2star} yields that the sequence
of $\R$-valued statistics 
$$
R_n((X_1,Y_1),\ldots,(X_n,Y_n)) 
 = 
\frac{1}{n} \sum_{i=1}^n 
 \Ls\bigl(X_i,Y_i,f_{\Ls,\D,\lb}(X_i)\bigr), \quad n\in\N,
$$
where 
$f_{\Ls,\D,\lb}
:=
\arg \min_{f\in H} \frac{1}{n} 
\sum_{j=1}^n  \Ls(X_j,Y_j,f(X_j)) + \lb \hhnorm{f}
$,
is uniformly qualitatively robust in a 
neighborhood $\mathcal{U}(\P_0)$ for every probability measure $\P_0 \in \PM(\cZ)$. 
Now we apply Theorem \ref{AC.Theorem3star}, which yields that
the sequence $\Law{n}{R;\P_n}$ of bootstrap SVM estimators of $\Law{n}{R;\P}$ is qualitatively robust for all $\P_0\in\PM(\cZ,\B(\cZ))$, which completes the proof.
\end{proofof}

\bibliography{//home//andreas//tex//christmann}

\end{document}